\documentclass[a4paper,12pt,leqno]{amsart}

\usepackage{amssymb,hyperref}
\usepackage[latin1]{inputenc}
\theoremstyle{plain}
\newtheorem{maintheo}{Theorem\!\!}

\newtheorem{theo}{Theorem}[section]

\newtheorem{cor}[theo]{Corollary}
\newtheorem{lem}[theo]{Lemma}

\theoremstyle{definition}

\theoremstyle{remark}
\newtheorem{rem}[theo]{Remark}

\numberwithin{equation}{section}

\newcommand{\thismonth}{\ifcase\month\or
  January\or February\or March\or April\or May\or June\or July\or
  August\or September\or October\or November\or December\fi
  \space\number\year}
\DeclareMathAlphabet{\mathrmsl}{OT1}{cmr}{m}{sl}

\newcommand{\oper}[3][n]{\newcommand{#2}{\mathop{\mathrm{#3}}%
\ifx n#1\nolimits\else\limits\fi} }
\newcommand{\rsoper}[3][n]{\newcommand{#2}{\mathop{\mathrmsl{#3}}%
\ifx n#1\nolimits\else\limits\fi} }
\newcommand{\RM}{{\mathbb R}}

\newcommand{\TM}{{\mathbb T}}

\newcommand{\vol}{\operatorname{vol}}
\newcommand{\Scal}{\operatorname{Scal}}
\newcommand{\Ric}{\operatorname{Ric}}
\newcommand{\tr}{\operatorname{tr}}

\renewcommand{\geq}{\geqslant}
\renewcommand{\leq}{\leqslant}
\newcommand{\cerc}{{\mathbb S}}

\renewcommand{\exp}{\operatorname{e}}

\newcommand{\proofof}[1]{\end{#1}\begin{proof}}

\newcounter{mnotecount}[section]
\renewcommand{\themnotecount}{\thesection.\arabic{mnotecount}}
\newcommand{\mnote}[1]
{\protect{\stepcounter{mnotecount}}$^{\mbox{\footnotesize  $
      \bullet$\themnotecount}}$ \marginpar{\raggedright\tiny\em
    $\!\!\!\!\!\!\,\bullet$\themnotecount: #1} }

\begin{document}

\title[Nonnegative Ricci curvature]{Conformally flat manifolds\\ with
nonnegative Ricci curvature}
\author{Gilles Carron}
\address{Laboratoire Jean Leray\\ UMR 6629 du CNRS\\ Universit\'e de Nantes\\
France}
\email{Gilles.Carron@math.sciences.univ-nantes.fr}
\author{Marc Herzlich}
\address{Institut de Math\'ematiques et Mod\'elisation de Montpellier\\ 
UMR 5149 du CNRS\\ Université Mont­pellier~II\\ France}
\email{herzlich@math.univ-montp2.fr}

\begin{abstract}
We show that complete conformally flat manifolds of dimension $n\geq 3$
with nonnegative Ricci curvature enjoy nice rigidity properties: they are 
either flat, or locally isometric to a product of a sphere and a line, or 
are globally conformally equivalent to $\RM^n$ or to a spherical spaceform
$\cerc^n/\Gamma$. This extends previous results due to Q.-M.~Cheng, M.-H.~Noronha,
B.-L. Chen and X.-P. Zhu, and S.~Zhu.
\end{abstract}
\thanks{Both authors are supported in part by an {\sc aci} program of the French 
Ministry of Research, and the {\sc edge} Research Training Network 
{\sc hprn-ct-2000-00101} of the European Union.}

\maketitle

\bigskip

In this note, we study complete conformally flat manifolds with nonnegative Ricci
curvature. It is well known in dimension $2$ that the sphere, the plane and
their quotients are the only surfaces that can be endowed with a metric of 
nonnegative curvature. As Riemannian surfaces are always conformally flat, it 
seems natural
to look at higher dimensional analogues of this fact. R.~Schoen and S.~T.~Yau 
showed in
\cite{kleinian} that conformal flatness together with nonnegative scalar curvature 
(or any variant of it involving the Yamabe constant) still allows much flexibility.
On the contrary, if one concentrates on stronger curvature conditions such as Ricci
curvature bounds, one might expect that they put some quite strong restrictions on 
the manifold. 

\smallskip

For instance, in case the manifold is closed, a characterization of the 
spaceforms and the quotients of $\cerc^1\times\cerc^{n-1}$ has been 
obtained by M.-H.~Noronha \cite{noronha-lcf}. If one adds extra assumptions like 
constant positive scalar curvature, further restrictions
can moreover be obtained, see for instance Q.-M. Cheng \cite{cheng-confplat} for 
results and references. 

\smallskip

More generally, and without any compactness assumption, S.~Zhu proved in 
\cite{zhu-lcf} that the {\sl universal covering} of a complete conformally flat
manifold of nonnegative Ricci curvature either is conformally equivalent to 
$\cerc^n$ or $\RM^n$, or is isometric to $\RM\times\cerc^{n-1}$. 
In the same vein, B.-L. Chen and X.-P. Zhu proved in \cite{chen-zhu} that 
the only complete non-compact conformally flat manifolds with nonnegative 
Ricci curvature and fast curvature decay at infinity are the complete  
flat manifolds. One motivation for this result is that it stands as an analog 
on real manifolds of well-known rigidity and gap phenomena on K\"ahler 
manifolds with nonnegative holomorphic bisectional curvature 
\cite{chen-zhu-kahler,mok-siu-yau,ni-vanishing,ni-mono,nishitam-poisson,ni-tam-plurisub}.

\medskip

The goal of this paper is to complete the above results by exhibiting a 
full classification of the possible geometries of complete conformally flat 
manifolds with nonnegative Ricci curvature. As we shall see, quite a bit more
can be said in this rather general setting on the topology of those manifolds, 
and also on their geometry. More precisely, we prove:

\smallskip

\begin{maintheo}
Let $(M,g)$ be a complete conformally flat manifold of dimension $n\geq 3$
with nonnegative Ricci curvature. Then, exactly one of the following holds:

{\upn 1.} $M$ is \emph{globally conformally equivalent} to $\RM^n$ with a 
conformal \emph{non-flat} metric with nonnegative Ricci curvature;

{\upn 2.} $M$ is \emph{globally conformally equivalent} to a spaceform of 
positive curvature, endowed with a conformal metric with nonnegative Ricci 
curvature;

{\upn 3.} $M$ is \emph{locally isometric} to the cylinder $\RM\times\cerc^{n-1}$;

{\upn 4.} $M$ is \emph{isometric} to a complete flat manifold.
\end{maintheo}

\medskip

Obviously, the last three cases do appear. In particular, since the round metric 
on spherical spaceforms has positive curvature, any small conformal deformation 
will keep nonnegative Ricci curvature as well, hence the second class is rather 
large. We shall show below that the first case also occurs, by exhibiting
explicit examples of non-flat globally conformally flat metrics with nonnegative
Ricci curvature on $\RM^n$. 

\medskip

The most interesting part of our theorem lies of course in the dichotomy between
the first two cases, where deformations are permitted, and the last two cases, 
where some strong rigidity is obtained. The philosophy of our result is then : 
a complete conformally flat manifold with nonnegative Ricci curvature is either 
(globally) {\sl topologically simple} (diffeomorphic to a vector space or a 
quotient of a sphere), or it is (locally) {\sl metrically rigid}.
Hence, only very general assumptions on the conformally flat manifold are enough 
to yield strong constraints either on its topology, or on its geometry. If one
compares with S.~Zhu's classification \cite{zhu-lcf} of the possible universal 
coverings, one sees that in the non-compact case ({\it i.e.} when quotients
of the sphere are excluded), the presence of a nontrivial fundamental group
implies metric rigidity.

\medskip

The vector space $\RM^n$ is then the only space that admits non-flat complete
conformally flat metrics with nonnegative Ricci curvature. As already alluded
to, B.-L.~Chen and X.-P.~Zhu \cite{chen-zhu} showed that, in case the metric
satisfies fast curvature decay assumptions at infinity, the metric is flat. It
is well known under very general circumstances that fast curvature decay at
infinity implies strong constraints on the topology at infinity, but our work 
implies further that, in the conformally flat case, either the metric is already
flat or the topology is precisely that of $\RM^n$. The examples we exhibit 
below of non-flat conformally flat metrics with nonnegative Ricci curvature on 
$\RM^n$ support the idea that the decay rate chosen by  B.-L.~Chen and X.-P.~Zhu 
\cite{chen-zhu} is indeed the optimal one.
It remains an interesting question to understand if further constraints 
can be deduced on non-flat conformally flat metrics with nonnegative Ricci 
curvature on $\RM^n$

\bigskip

\begin{small}
{\flushleft\it Acknowledgements}. The authors are grateful to J\'er\^ome~Droniou 
and Jacques~Lafontaine for useful discussions. Many thanks are due to Jeff 
Viaclovsky who pointed out S.~Zhu's paper to us. 
\end{small}

\bigskip

\section{Non trivial conformally flat metrics on $\RM^n$\\ with nonnegative Ricci
curvature}

\smallskip

In this short section, we show that the Theorem above is indeed
the best one can hope for. As a matter of fact, 
we exhibit examples of (rotationnally symmetric) metrics on $\RM^n$
that are {\sl non flat}, but {\sl complete, conformally flat and with nonnegative 
Ricci curvature}. Hence this may exist, whereas existence of any analogous metric
on quotients is forbidden by our result. A by-product is the optimality of the 
curvature decay rate imposed by B.-L.~Chen and X.-P.~Zhu 
in \cite{chen-zhu} to get a flat metric: our
example lies exactly on the threshold where the result in \cite{chen-zhu} does not 
hold anymore. 

{\flushleft\sl Construction of the example}. 
Let $f$ be a real function. From Besse's book \cite[formula 1.159, page 59]{besse} 
we get the expression for the Ricci curvature of the metric
$\exp^{2f}\!g_0$ ($g_0$ the euclidean metric) is 
\begin{equation}
\Ric = - (n-2) (Ddf - df\otimes df) + (\Delta f - (n-2)|df|^2 ) g_0 
\end{equation}
which can be rewritten, in case $f=f(r)$ is a radial function and $h$ is
the unit round metric on $\cerc^{n-1}$, as
\begin{equation}
\begin{split}
\Ric = & -(n-2) \left( f'' dr^2 - (f')^2 dr^2 + r f' h \right) \\
       & \ \ + \left( - f'' - (n-1)r^{-1} f' -(n-2) (f')^2 \right)\, (dr^2 
+ r^2 h) \\
     = & - (n-1) \left( f'' + r^{-1} f' \right)\, dr^2  \\
& \ \ - \left( f'' + (2n-3) r^{-1} f' + (n-2) (f')^2 \right)\, r^2 h .
\end{split}
\end{equation}
The Ricci tensor is then nonnegative iff.
\begin{equation}\label{eq1}
f'' + r^{-1}f' \leq 0
\end{equation}
and
\begin{equation}\label{eq2}
f'' + (2n-3)r^{-1}f' + (n-2) (f')^2 \leq 0 .
\end{equation}
It will be convenient to work with $f'(r) = r^{-1}a(r)$. With this choice, the
conditions for Ricci to be nonnegative become:
\begin{equation}\label{neweq1}
a' \leq 0
\end{equation}
and 
\begin{equation}\label{neweq2}
a' + (n-2) r^{-1} (2a + a^2) \leq 0 .
\end{equation}
A function satisfying (\ref{neweq1}--\ref{neweq2}) is easily
found: for instance, one can choose any smooth function $a$
such that
\begin{equation}
a(r) = 0  \ \ \forall r \in [0,1] \ , \ 
a(r) = -\frac{1}{2} \ \ \forall r \geq 2
\end{equation}
so that $a'(r)\leq 0$ and $-2 \leq a(r) \leq 0$ everywhere. We then
let 
\begin{equation}
f(r) = \int_0^r \frac{a(s)}{s} \, ds
\end{equation}
Hence the metric $\exp^{2f}g_0$ is conformally flat and has nonnegative
Ricci curvature. This metric behaves as $r^{-1}\, (dr^2 + r^2 h)$ at
infinity, it is then complete since
\begin{equation}
\int_1^{\infty} \frac{dr}{\sqrt{r}} = + \infty .
\end{equation}

\medskip

{\flushleft\sl Optimality of the decay condition of B.-L.~Chen and X.-P.~Zhu}
\cite{chen-zhu}. 
Let us analyse a little further the examples above. It is clear
that one may take:
\[
\lim_{r\rightarrow\infty} a (r) = - (1 - \varepsilon)\ \ \textrm{ for any } \ 
\varepsilon > 0,
\]
and the metric one obtains is still complete. A short computation leads to
an estimate of the Ricci curvature:
\[
\Ric \simeq  (1- \varepsilon^2) r^{-2\epsilon}.
\]
But geodesic distance from the origin in this metric is $s= r^{-\epsilon}$ around 
infinity, hence
\[
|\Ric| = O(s^{-2}).
\]
In \cite{chen-zhu}, B.-L.~Chen and X.-P.~Zhu proved that a 
complete non-compact conformally
flat manifold with nonnegative Ricci curvature is necessarily flat if the scalar
curvature is bounded, and if, for geodesic balls centered at some fixed origin,
\[ \frac{s^2}{\vol B(s)} \int_{B(s)} |\Ric|  \longrightarrow 0 \ \ \textrm{ as }
s \textrm{ goes to infinity}. \]
In particular, this is satisfied if the following simplest conditions hold
true: 
\[ \vol B(s)\geq C s^n, \ \ s^{2+\delta} |\Ric | \longrightarrow 0 \ \ \textrm{ as }
s \textrm{ goes to infinity, for some } \delta > 0 . \]
As a result, we see that the assumptions in B.-L.~Chen and X.-P.~Zhu \cite{chen-zhu}
is close to be optimal on manifolds conformally equivalent to $\RM^n$.

\bigskip

\section{The proof}

\smallskip

The proof of the Theorem is divided into two parts: the first one is a 
classification of the possible holonomy coverings of $(M,g)$ be a complete 
(compact or not) conformally flat manifold with nonnegative Ricci curvature.
This is basically S.~Zhu's result \cite{zhu-lcf}, 
but for sake of completeness, and also because
this is quite a short proof, we have provided below the arguments leading to 
the results we will precisely need. Our arguments are roughly similar to those 
of M.-H.~Noronha \cite{noronha-lcf} and S.~Zhu \cite{zhu-lcf},
although they differ at some places. It then turns out that the classification of
the holonomy coverings is a classification of the universal coverings, and
it then remains to show that, in the presence of a fundamental group, one
gets metric rigidity (unless we are in the spherical spaceform case).

\smallskip

The very first tool one needs to achieve the holonomy covering classification
is R.~Schoen and S.~T.~Yau's 
analysis \cite{kleinian} of conformally flat manifolds with 
nonnegative scalar curvature. 
This applies in dimensions $n\geq 7$ without any further restrictions, 
and to manifolds of dimensions $3\leq n\leq 6$ provided
that a very general positive mass holds true for manifolds with one asymptotically 
flat end and a huge (possibly infinite) number of other complete ends. 

This last result is considered as highly probable but has unfortunately remained 
unpublished so far (see \cite[\S 4]{kleinian} and the appendix for details on 
this subtle point). Luckily enough, we only need Schoen-Yau's analysis to hold 
true in our setting where extra conditions on the Ricci curvature are available,
and it turns out that a direct proof can be obtained using this extra assumption. 
For sake of 
completeness, we have included a description of the elements of that proof in an
appendix to this paper. Note also that U.~Christ and J.~Lohkamp recently announced 
a proof of a general positive mass theorem with new techniques that might lead to 
the desired result \cite{christ-lohkamp-annonce}.

\smallskip
    
We now come back to our main proof: 
Schoen-Yau's analysis then implies that the holonomy covering of $M$ 
is a domain $\Omega$ included in the complement $\Omega(\Gamma)$ 
of the limit set of the fundamental group $\Gamma = \pi_1(M)$ of $M$ in 
the sphere $\cerc^n$.

\smallskip
 
We now look at $\Omega$ : it is endowed with a complete metric conformal to 
the round metric and, due to nonnegativity of Ricci curvature, it has 
sub-euclidean volume growth:
\[ \vol\, B(r) \leq C r^n .\]
One can now apply \cite[Prop. 2.2]{carron-mh}, which shows that the
complement of $\Omega$ in the sphere is of $n$-capacity zero, hence
of zero Hausdorff dimension. Moreover, manifolds with nonnegative Ricci 
curvature have a finite number of ends \cite{liu-ric}, and, arguing as
in \cite[Proof of Theorem 2.1]{carron-mh}, we can conclude that $\Omega$
must be the complement of a finite number of points in the sphere $\cerc^n$. 
Three cases can now occur: 

\smallskip

{\flushleft (1)} First of all, $\Omega$ might be the whole sphere. Then
$M$ is compact, with a metric $g$ conformal to the round metric and with 
nonnegative Ricci curvature. As already noticed, 
there are plenty of such metrics, such as,
{\sl e.g.}, $C^2$-small conformal deformations of the round metric. 

\smallskip

{\flushleft (2)} Assume now that $\Omega$ has at least two ends. Since it has
nonnegative Ricci curvature, Cheeger-Gromoll theorem \cite{CG1} 
implies that $\Omega$ splits as a Riemannian product $\RM\times Y$. Remarks due
to Lafontaine \cite{jaclaf-confplat} assert that a conformally flat manifold
that is a Riemannian product with a line is necessarily a Riemannian product of a
constant curvature manifold with the line. Hence, $Y$ must be of constant 
curvature, and as $\Omega$ is simply connected, $Y=\RM^{n-1}$ or $\cerc^{n-1}$ 
(hyperbolic space would lead to negative Ricci curvature). Since $\Omega$ has 
at least two ends, we end up with $Y=\cerc^{n-1}$ as the sole possibility in
this case.

\smallskip

{\flushleft (3)} If $\Omega$ has only one end, it forces it to be diffeomorphic
to $\RM^n$, endowed with a metric conformal to the euclidean metric and with 
nonnegative Ricci curvature. Either $M$ itself is diffeomorphic to $\RM^n$ 
(and we have seen in the previous section that the metric is not necessarily
flat) or $M$ is different from its universal covering $\Omega$. From now on, 
we will assume that we are in that last case, and we will show that $g$ is
flat.

\smallskip

If $\Omega(\Gamma) = \cerc^n$, it implies that $\Gamma$ is a group of conformal
diffeomorphisms of the sphere acting properly discontinuously without fixed
points on the whole sphere. It is then a finite group, and it must fix the 
missing point $\Omega(\Gamma) - \Omega$. Any conformal map of finite order has 
at least two fixed points (if not they would be parabolic, and this would imply 
both $\Omega(\Gamma) =\Omega$ and infinite order) and $\Gamma$ cannot act without 
fixed points on $\Omega$, hence this case is excluded.

\smallskip
 
We are left with the case where $\Omega(\Gamma) = \Omega$ is diffeomorphic to
$\RM^n$ and $\Gamma$ acts properly discontinuously without fixed points as 
conformal diffeomorphisms of the flat space. But the only possibility is then 
that $\Gamma$ is a subgroup of the group of Euclidean isometries of flat space.

\smallskip

Let us note now that if $M$ itself is compact, then it is (isometrically) a 
compact flat manifold. Indeed, $M$ is diffeomorphic to a flat manifold, and
is endowed with a conformal metric $g$ with nonnegative Ricci curvature.
By Bieberbach theorem \cite{wolf} it is covered by a torus, and the Bochner 
technique \cite[4.37]{ghl3ed} (or, alternatively, the Gromov-Lawson scalar 
curvature 
obstruction \cite{gromov-lawson-ihes}) applies to show that $g$ itself is flat.

\smallskip

We now consider the case where $M$ is non-compact. 
By Bieberbach theorem again,  there exists a finite cover $\widehat{M}$
of $M$, which is conformally equivalent to a product of a flat torus and a flat 
space : 
\[ \widehat{M} = \TM^k \times \RM^{n-k}. \]
Moreover, $M$ is obtained from $\widehat{M}$ by a finite group acting isometrically
on the flat metric. Top conclude at this point, one could try to look at an 
extension of Gromov-Lawson's obstruction to nonnegative scalar curvature for
general non-compact Riemannian manifolds of type $\TM^k \times \RM^{n-k}$.
However, a careful reading of \cite[pp. 95--96]{gromov-lawson-ihes} shows that it 
is useless to hope for such an extension unless $k=n-1$. We shall then need
an {\sl ad hoc} proof in the conformally flat case.

Let us now study the different possible values for $k$ (with $1\leq k\leq n-1$ 
as the case $k=n$ has already been treated above and the case $k=0$ bears no
restriction, as we know from section~1).

\smallskip

{\flushleft (1)} If $k=n-1$, then $\widehat{M}$ has two ends, and Cheeger-Gromoll's
theorem \cite{CG1} applied again together with Lafontaine results 
\cite{jaclaf-confplat} shows that $M$ is isometric to a flat 
spaceform.

\smallskip

{\flushleft (2)} Suppose now that $1\leq k\leq n-2$. 
Using the computations done above, the Ricci curvature
of the metric $g=\phi^{-2}g_0$ ($g_0$ the euclidean metric on 
$\TM^k\times\RM^{n-k}$) is
\[ \Ric^{g} = (n-2) \frac{Dd\phi}{\phi} - \left( \frac{\Delta\phi}{\phi}
+ (n-1)\frac{|d\phi|^2}{\phi^2} \right) g_0 .\]
Suppose now that $\widehat{M}= \TM^k\times \RM^{n-k}$ for $0 < k < n-1$, and
use polar coordinates on the $\RM^{n-k}$-factor to get
\begin{equation*}\begin{split}
\phi\Ric^{g}(\partial_r,\partial_r) \ = \  & (n-2)(\partial_r^2\phi) \\ 
& - \left( - (\partial_r^2\phi) - \frac{n-k-1}{r}\partial_r\phi
+ (n-1)\frac{|d\phi|^2}{\phi} + \frac{1}{r^2}\Delta_{\cerc^{n-k-1}}\phi
+ \Delta_{\TM^k}\phi \right)
\end{split}\end{equation*}
so that 
\begin{equation}\label{eq:1} 
0 \leq (n-1)(\partial_r^2\phi) + \frac{n-k-1}{r}\partial_r\phi - 
(n-1)\frac{(\partial_r\phi)^2}{\phi} - \frac{1}{r^2}\Delta_{\cerc^{n-k-1}}\phi
- \Delta_{\TM^k}\phi .
\end{equation}
Similarly,
\begin{equation*}\begin{split}
\phi\,\tr_{\TM^k}\Ric^{g} \ = \ & -(n-2)\Delta_{\TM^k}\phi \\ 
& - k\left( - (\partial_r^2\phi) - \frac{n-k-1}{r}\partial_r\phi
+ (n-1)\frac{|d\phi|^2}{\phi} + \frac{1}{r^2}\Delta_{\cerc^{n-k-1}}\phi
+ \Delta_{\TM^k}\phi \right),
\end{split}\end{equation*}
so that
\begin{equation}\label{eq:2}
0 \leq (\partial_r^2\phi) + \frac{n-k-1}{r}\partial_r\phi - 
(n-1)\frac{(\partial_r\phi)^2}{\phi} - \frac{1}{r^2}\Delta_{\cerc^{n-k-1}}\phi
- (1+\frac{n-2}{k})\,\Delta_{\TM^k}\phi .
\end{equation}
We now integrate equations (\ref{eq:1}) and (\ref{eq:2}) on 
$\cerc^{n-k-1}\times\TM^k$, and we let 
\[ \bar\phi = \int_{\cerc^{n-k-1}\times\TM^k} \phi .\]
This yields
\begin{equation}\label{eq:12'}\begin{split}
0 \ \leq & \ (n-1)\,\partial_r^2\bar\phi + \frac{n-k-1}{r}\partial_r\bar\phi - (n-1)
\,\int_{\cerc^{n-k-1}\times\TM^k} \frac{(\partial_r\phi)^2}{\phi} \ ; \\
0 \ \leq & \ \partial_r^2\bar\phi + \frac{n-k-1}{r}\partial_r\bar\phi - (n-1)
\,\int_{\cerc^{n-k-1}\times\TM^k} \frac{(\partial_r\phi)^2}{\phi} \ .
\end{split}\end{equation}
By H\"older inequality,
\begin{equation*}\begin{split}
\left(\partial_r\bar\phi \right)^2 & = \left( \int_{\cerc^{n-k-1}\times\TM^k}
\partial_r \phi \right)^2 
= \left( \int_{\cerc^{n-k-1}\times\TM^k} \phi^{1/2}\phi^{-1/2}\partial_r \phi 
\right)^2 \\
 & \leq \left( \int_{\cerc^{n-k-1}\times\TM^k} \phi\right)
\, \left( \int_{\cerc^{n-k-1}\times\TM^k} \frac{(\partial_r \phi)^2}{\phi} 
\right)
 = \bar\phi \, \left( \int_{\cerc^{n-k-1}\times\TM^k} 
\frac{(\partial_r \phi)^2}{\phi} \right). 
\end{split}\end{equation*}
Hence our pair of equations are transformed into the following:
\begin{equation}\label{eq:12''}\begin{split}
0 \ \ \leq & \ (n-1)\,\partial_r^2\bar\phi + \frac{n-k-1}{r}\partial_r\bar\phi - 
(n-1)\,\frac{(\partial_r\bar\phi)^2}{\bar\phi} \ ; \\
0 \ \ \leq & \ \partial_r^2\bar\phi + \frac{n-k-1}{r}\partial_r\bar\phi - (n-1)
\,\frac{(\partial_r\bar\phi)^2}{\bar\phi} \ .
\end{split}\end{equation}
Using this and the change of variables
\[ r = \exp^{t} , \ \ {\textit i.e.} \ r\partial_r = \partial_t, \ 
r^2\partial_r^2 = \partial_t^2 - \partial_t ,\]
and denoting $\psi(t) = \bar\phi(\exp^t)$ for $t\in\RM$, one gets:
\begin{equation}\label{eq:1*}
(n-1)\, \partial_t^2\psi - k\, \partial_t\psi 
- (n-1)\,\frac{(\partial_t\psi)^2}{\psi} \ \geq \ 0 \ , 
\end{equation}
and
\begin{equation}\label{eq:2*}
\partial_t^2\psi + (n-k-2)\, \partial_t\psi 
- (n-1)\,\frac{(\partial_t\psi)^2}{\psi} \ \geq \ 0\ . 
\end{equation} 
From equation (\ref{eq:1*}) one gets, if $k>0$:
\begin{equation}\label{eq:concl}
 \partial_t\psi \ \leq \ \frac{n-1}{k} \, \partial_t^2\psi - \frac{n-1}{k} 
\frac{(\partial_t\psi)^2}{\psi} . 
\end{equation}
Notice now that $k\leq n-2$, so that injecting (\ref{eq:concl}) into equation 
(\ref{eq:2*}) leads to:
\begin{equation}\label{eq:3} 
\partial_t^2\psi - \left( 1 + \frac{k(n-2)}{k+(n-k-2)(n-1)} \right)
\frac{(\partial_t\psi)^2}{\psi} \ \geq \ 0. 
\end{equation}
Define now
\begin{equation}\label{eq:alpha} 
\alpha = - \, \frac{k(n-2)}{k+(n-1)(n-k-2)} \ < \ 0 ,\end{equation}
so that
\begin{equation*}\begin{split}
\partial_t^2(\psi^{\alpha}) & = \partial_t \left( \alpha\, 
\psi^{\alpha-1}\,\partial_t\psi \right) \\
 &  = \alpha \left( \psi^{\alpha -1}\partial_t^2\psi + (\alpha-1)\,\psi^{\alpha-2}
\,(\partial_t\psi)^2 \right) \\
 & = \alpha\,\psi^{\alpha-1}\left( \partial_t^2\psi - 
(1-\alpha)\,\frac{(\partial_t\psi)^2}{\psi} \right).
\end{split}\end{equation*}
As $\alpha<0$ and $\psi>0$, this implies that $\psi^{\alpha}$ is a positive 
concave function on $\RM$. But constants are the only positive concave functions 
on $\RM$ that do not
change signs, so that the function $\bar\phi$ is necessarily a constant in $r$. 

\medskip

We now define 
\[ \Psi : x \in \RM^{n-k} \longmapsto \int_{\{ x\}\times\TM^k} \phi \ \ .\] 
Using $x$ in $\RM^{n-k}$ as the origin, the work done above 
shows that, for each $r>0$,
\begin{equation*}\begin{split} \Psi (x) \ & = 
\ \lim_{s\rightarrow 0}\, \frac{1}{\vol\,\cerc^{n-k-1}}\, \bar\phi(s) \\
& = \ \frac{1}{\vol\,\cerc^{n-k-1}}\, \bar\phi (r) \\
& = \frac{1}{\vol\,\cerc^{n-k-1}}\, \int_{\{ r\}\times\cerc^{n-k-1}\times\TM^k} 
\phi \\
& = \ \frac{1}{\vol S_x(r)} \,
\int_{S_x(r)} \Psi \ \ ,\end{split}\end{equation*}
or, equivalently: the function $\Psi$ is harmonic on $\RM^{n-k}$.

\medskip

We are now in a position to end the proof. The condition for the scalar curvature 
to be nonnegative is :
\begin{equation}\label{eq:scal} 
0 \ \leq \ \Scal^g \ = \ - 2(n-1)\Delta\phi - n(n-1)\frac{|d\phi|^2}{\phi} ,
\end{equation}
so that 
\begin{equation}\label{eq:harmonic} \Delta\phi \ = \ \Delta_{\RM^{n-k}} \phi 
\ + \ \Delta_{\TM^k}\phi \ \leq \ 0 \ .\end{equation}
Integrating on  $\TM^k$, the first term in (\ref{eq:harmonic})
vanishes since $\Psi$ is harmonic,
whereas the second vanishes by definition of the Laplacian on $\TM^k$.
As a result, $\phi$ itself must be harmonic and the sign of the 
scalar curvature equation (\ref{eq:scal}) shows that $d\phi$ must vanish.
Hence $g$ must be flat. \qed

\bigskip

\begin{rem}
Notice that the end of the proof given above works exactly when we expect the 
result to be true. Indeed, it fails if $\alpha$ in equation (\ref{eq:alpha}) 
vanishes, {\sl i.e.} if  $n=2$ (surface case, where an alternative proof must be 
given) or if $k=0$ (no torus factor, where the metric might be non-flat, in 
agreement with the example given in the first section).
\end{rem} 

\medskip

\begin{rem}
For some values of $k$ and $n$ (namely $k\geq n-2$), a more geometric proof 
can be obtained along the following lines. One has
\[ M = \left( \TM^k \times \RM^{n-k} \right) / \Delta  \]
for some finite group $\Delta$, so that $\Gamma = \pi_1(M)$ has
polynomial growth of order $k$.
Moreover, from Bishop-Gromov classical comparison theorem, volume growth
is polynomial of order $n$ at most. It was proved by Anderson 
\cite{anderson-riccipositif} that if $\widehat{M}=\TM^k \times \RM^{n-k}$ 
has volume growth 
\[ \vol (B(r)) \geq C r^{p} \]
then the group $\widehat{\Gamma} =  \Gamma/\Delta$ has
polynomial growth $n-p$ at most, hence $\Gamma$ itself. 
This implies $p\leq n-k$. 
Applying Anderson's result once more, polynomial growth $k$ of
$\widehat{\Gamma}$ and volume growth 
\[ \vol (B(r)) \leq C' r^{n} \]
on $\Omega$ implies that $\widehat{M}$ has volume growth
\[  \vol (B(r)) \leq C' r^{n-k}.  \] 
Now, if the curvature is bounded, and if $n-k\leq 2$, Anderson has proved
in the same paper using minimal surfaces that $\widehat{M}$ must split off 
a factor $\RM$ at least. Lafontaine's result \cite{jaclaf-confplat} then 
ends this remark. Unfortunately, there seem to be no way to extend this 
conceptual argument to higher codimensions.
\end{rem}

\bigskip

\section*{Appendix}

\medskip

In the last proof, we have used the following result of R. Schoen and S.T. Yau:

\smallskip

\begin{theo}[Schoen-Yau \cite{kleinian}, Theorems 4.1--4.5]\label{SY}
Suppose $(M,g)$ is complete, locally conformally flat with nonnegative scalar
curvature, then the developing map :
$$\Phi : \tilde M\rightarrow \cerc^n$$ is injective, the fundamental group
$\pi_1(M)$ is isomorphic to
a discrete subgroup of $C_n$ \upn{(}the conformal group of $\cerc^n$\upn{)}, 
and $\Phi(\tilde M)$ is a domain $\Omega$ in the complement $\Omega(\Gamma)$ 
of the limit set of the action of $\pi_1(M)$ on the sphere $\cerc^n$.
\end{theo}

\smallskip

The proof of this result is done in \cite{kleinian} in
two steps: a first argument covers the case 
$n\ge 7$; another argument based on the positive mass theorem covers the case 
$3\le n\le 6$. However, the authors remark that "{\it for this application it is 
necessary to extend the positive mass theorems to the case of complete manifolds; 
that is assuming that the manifold has an asymptotically flat end and other ends 
which are merely complete. This extension will be carried out in a future work}"
\cite[p.~65]{kleinian}. 
It is widely believed in the Riemannian Geometry community that the needed 
extension of the positive mass theorem is true, but, unfortunately, no such
generalisation has been published so far. Of course,
when the manifold is assumed to be spin, the Witten proof of the positive mass 
theorem \cite{parker,parker-taubes,witten}
extends to this general framework but the spin assumption plays an essential
role there. Note however that every orientable $3$-manifold is spin. Hence, 
Schoen-Yau's analysis of conformally flat manifolds and Kleinian groups
\cite{kleinian} for sure holds in dimension $3$ as we only need to 
apply the positive mass theorem to the universal covering (mass increases as 
one takes quotients).

\medskip

To avoid in our setting the problems caused by the absence of any written proof
of the relevant positive 
mass theorem, it then remains to prove the adequate result for conformally 
flat manifolds with nonnegative Ricci curvature and dimensions $4\le n\le 6$. 

\smallskip

According to the Bishop-Gromov inequality, we know that on any covering of such a
manifold the volume of geodesic balls is less than the volume of the euclidean balls
of the same radius.  We will show in this appendix that the argument given  in 
\cite{kleinian} by Schoen and Yau to prove their theorem in dimension larger than 
$7$ works in dimensions $4$, $5$, and $6$ under just that very assumption 
on volume growth that
can be obtained from nonnegativity of the Ricci curvature. More precisely, we
prove below the following version of Theorem 3.1 in \cite{kleinian}:

\smallskip

\begin{theo}[Schoen-Yau's Theorem 3.1 revisited] Assume that $(M,g)$ is a complete 
Riemannian conformally flat manifold of dimension $n\in [4,6]$ such that
$\Scal_g\ge0$, and, for some $x\in M$, there exist $C$ such that
$\vol B(x,r)\le C r^n$ for every $r>0$. 

Then if $\Phi\,:\, M\rightarrow \cerc^n$ is a conformal immersion then $\Phi$ is
injective and yields a conformal diffeomorphism of $M$ onto $\phi(M)$. Moreover
the boundary of $\Phi(M)$ has zero Newtonian capacity.
\end{theo}

\smallskip

Our proof closely follows \cite{kleinian}, differing only at the points where 
dimension arguments enter the picture. We then refer to \cite[p.~58--61]{kleinian} 
for the overall strategy, and focus below only on the steps where some modification 
is needed. In particular, we refer to \cite{kleinian} for the proof of the fact
that Theorem 2.4 leads to the following final result:

\smallskip

\begin{cor}
Suppose $(M,g)$ is complete, locally conformally flat with nonnegative Ricci
curvature then the developing map :
$$\Phi\, :\,  \tilde M\rightarrow \cerc^n$$ is injective, $\pi_1(M)$ is 
isomorphic to a discrete subgroup of $C_n$ \upn{(}the conformal group of 
$\cerc^n$\upn{)}, and $\Phi(M)$ is a domain $\Omega$ included in the complement 
$\Omega(\Gamma)$ of the limit set of $\pi_1(M)$ in $\cerc^n$.
\end{cor}

\smallskip

{\flushleft\it Proof of Theorem 2.4. --} We will first show how the original proof 
of Schoen and Yau can be easily adapted for the case $n\ge 5$; then we examine the 
case $n=4$ which needs a more elaborate argument.

\smallskip

The core of the proof of Theorem 3.1 in \cite{kleinian} lies in the following
construction: let $o\in M$ be a point and $G_o$ be the Green kernel of the 
conformal Laplacian of $(M,g)$ with pole at $o$, and $\bar G_o$ be the pulled-back  
(by $\Phi$) of the Green kernel of the conformal Laplacian on the sphere with pole 
at $\Phi(o)$, normalized so that $\bar G_o-G_o$ is smooth near $o$. 
We let $v=G_o/\bar G_o$. The metric $\bar g =G_o^{4/(n-2)} g$ is flat (and 
incomplete) and the function $v$ is $\bar g$-harmonic and one has $0<v\le 1$. 
Then Schoen and Yau prove that their Theorem is obtained if one is able to
show that $v=1$ everywhere. It is that careful study of $v$, where the positive
mass theorem plays a role in small dimensions, that we shall replace.

\smallskip

When $n\ge 5$ and $\alpha=2(n-2)/n$, R. Schoen and S.T. Yau show in 
\cite[equation (3.4), p.~59]{kleinian} that for every $\sigma$ large enough :
$$\int_{B(o,\sigma/2)} |\bar\nabla v|^{\alpha-2}  |\bar\nabla  |\bar\nabla v||^2
dvol_{\bar g}\le \frac{C}{\sigma^2}\int_{B(o,3\sigma)\setminus B(o,\sigma/3)}
G_o^\alpha \ d\!\vol_g \ .$$
Moreover, they have also shown that for every $p>n/(n-2)$ 
\cite[Corollary 2.3]{kleinian}:
$$ \int_{M \setminus B(o,1)} G_o^p \, d\!\vol_g <  \infty \ ;$$
in particular as $n>4$ we have 
$$ \int_{M \setminus B(o,1)} G_o^2 \, d\!\vol_g <  \infty\  .$$
Using H\"older inequality and our assumption on volume growth, we then have
$$\int_{B(o,3\sigma)\setminus B(o,\sigma/3)}
G_o^\alpha \, d\!\vol_g \ \le 
\ C \sigma^2 \,\left(\int_{B(o,3\sigma)\setminus B(o,\sigma/3)}
G_o^2 \, d\!\vol_g\right)^{1-2/n} .$$
Hence, letting $\sigma\to \infty$, we obtain that $v$ is the constant function
equal to $1$. As already said, this is enough to continue the proof with the 
argument of \cite{kleinian}.

\medskip

It remains to consider the case of the dimension $4$ where, unfortunately, 
$\alpha=1$. The first results we need are :

\smallskip

\begin{lem}\label{lem:first}
For $R>0$ large enough and for all $\varepsilon>0$,
\[ \int_{M\setminus  B(o,R)}\frac{ |\nabla G_o|^2}{G_o|\log G_o|^{1+\varepsilon}}
d\!\vol_g  \ + 
\ \int_{M\setminus  B(o,R)}\frac{\Scal_g G_o }{|\log G_o|^{\varepsilon}}
\vol_g <\infty \]
\end{lem}

\smallskip

\begin{proof} This is inspired by \cite[proof of Proposition 2.4(iii)]{kleinian}.

Let $M=\cup_i U_i$ an exhaustion of $M$ by a sequence of nested bounded open subsets
with smooth boundaries, such that $o\in U_0$. Let $\varphi\in C^\infty(M)$ such
that $\varphi=1$ outside $B(o,R)$ and the support of $\varphi$ does not meet
$B(o,R/2)$ and let $G_o^i$ be the Green kernel of the confomal Laplacian on $U_i$
(with the Dirichlet boundary condition) of pole $o$.

We know that $G_o^i- G_o\to 0$ uniformly on compact subsets of $M$ and that 
$G_o^i\le G_o$. Moreover, nonnegativity of the scalar curvature implies that 
\begin{equation}\label{eq:scacal}
\Delta G_o^i+\Scal_g G^i_o/6\le 0 .
\end{equation}
Hence we can integrate by parts the quantity 
\[ \int_{U_i} \varphi  \langle 
\nabla G_o^i,\nabla(|\log G_o^i|^{-\varepsilon})\rangle 
\, d\!\vol_g
\] 
with $i$ large enough, so that $B(o,R) \subset U_i$. Note that it is shown in
Lemma 3.2 in \cite{kleinian} that 
\[ \lim_{x\rightarrow\infty} G_o(x) = 0 ,\]
so that one may take $i$ and $R$ large enough to ensure that $G_o^i\leq 1/2$
on the support of $\varphi$
Using inequality 
(\ref{eq:scacal}), and letting $i\to\infty$, this leads to :
\begin{equation*}\begin{split}
\varepsilon\,\int_{U_i\setminus  B(o,R)} & 
\frac{ |\nabla G^i_o|^2}{G^i_o|\log G^i_o|^{1+\varepsilon}}\, d\!\vol_g \\
& + \ \frac{1}{6} 
\int_{U_i\setminus B(o,R)}
\frac{\Scal_g G_o }{|\log G_o|^{\varepsilon}} \, d\!\vol_g
\le C \int_{B(o,R)\setminus B(o,R/2)}
|\nabla G^i_o||\log G^i_o|^{-\varepsilon}d\!\vol_g .
\end{split}\end{equation*}
The result follows from letting $i\to \infty$, because the right-hand side
is an integral on a domain that does not depend on $i$ and the convergence
of $G_o^i$ towards $G_o$ is in $C^{\infty}$ on compact subsets; 
moreover this is true for every $\varepsilon >0$.
\end{proof}

\smallskip

\begin{lem}\label{finite2} For any $R>0$,
\[ \int_{M\setminus  B(o,R)} \frac{G_o ^2}{|\log G_o|^{4}} d\!\vol_g<\infty .\]
\end{lem}

\smallskip

\begin{proof}
We notice that R. Schoen and S.T. Yau showed that the Yamabe invariant 
of $M$ is positive (in fact it is exactly that of the standard sphere, see 
Proposition 2.2 in \cite{kleinian}),
hence we have a constant $c>0$ such that for any $f\in C^\infty_0(M)$,
$$c \left( \int_M f^4 \, d\!\vol_g\right)^{1/2}\le 
\int_M [ |df|^2+\frac16 \Scal_g f^2 ]\, d\!\vol_g$$
In particular, we can apply this to $f = \varphi\sqrt{G_o^i}/\log G_o^i$ (with
$\varphi$ a cut-off function as above) and Lemma 
\ref{lem:first} with adequate choices of $\varepsilon$ shows the desired result. 
\end{proof}

\smallskip

Now we can adapt the argument of \cite{kleinian} to the case $n=4$:
we start again from inequality (3.4) of that paper:
$$\forall \phi\in C^\infty_0(M),\ 
\int_M \phi^2 |\bar\nabla v|^{-1}  |\bar\nabla  |\bar\nabla v||^2\,
dvol_{\bar g} \ \le\ \int_M  |d\phi|^2G_o |dv| \, d\!\vol_g .$$
Let $k$ an integer large enough; we take the following as a test function :
$$\phi_k(x)=\left\{
\begin{array}{lll}
1& {\rm if} & x\in B(o,2^k)\\
\frac{ \log\left(2^{2k}/r(x)\right)}{\log (2^k)} \ \leq 1 
& {\rm if} & x\in B(o,2^{2k})\setminus B(o,2^k)\\
0&{\rm if} & x\not \in B(o,2^{2k})\\
\end{array}
\right.$$
where $r(x)={\rm dist}(x,o)$.
As we have 
\[ 
G_o |dv|\le |dG_o|+G_o |d\log \bar G_o|,
\] 
we obtain :
\begin{equation}\label{eq:1and2}
\int_{B(o,2^k)}  |\bar\nabla v|^{-1}  |\bar\nabla  |\bar\nabla v||^2
dvol_{\bar g}\ \le\ \frac{c}{k^2} \int_{B(o,2^{2k})\setminus B(o,2^k)} [ |dG_o|+G_o
|d\log \bar G_o|\, ] \, \frac{d\!\vol_g}{r^2(x)}.
\end{equation}

\smallskip

We first treat the first term in the right-hand side of the previous estimate.
One has:
$$\int_{B(o,2^{2k})\setminus B(o,2^k)}  |dG_o| \frac{d\!\vol_g}{r^2(x)}\le
\left(\int_{B(o,2^{2k})\setminus B(o,2^k)}  |dG_o|^2d\!\vol_g\right)^{1/2}
\left(\int_{B(o,2^{2k})\setminus B(o,2^k)} 
\frac{d\!\vol_g}{r^4(x)}\right)^{1/2}.$$ 
Now the same proof as that of Lemma \ref{lem:first} above (see also 
\cite[Proposition 2.4(iii-a)]{kleinian}) shows that 
$$\int_{M\setminus B(o,1)} |dG_o|^2 \, d\!\vol_g<\infty \ , $$ 
and the upper bound on the volume growth yields 
$$\int_{B(o,2^{2k})\setminus B(o,2^k)} \frac{d\!\vol_g}{r^4(x)}\ \le \ C \, k$$ 
for a constant $C$ independent of $k$. Hence we get that the first term goes to 
zero when $k$ goes to infinity.

\smallskip

Now for the second term in the right-hand side of (\ref{eq:1and2}):
\begin{equation*}
\begin{split}
& \int_{B(o,2^{2k})\setminus B(o,2^k)} G_o |d\log \bar G_o| 
\frac{d\!\vol_g}{r^2(x)}  \\
& \ \ \leq \ \left( \int_{B(o,2^{2k})\setminus B(o,2^k)} G_o|\log G_o|^{2} 
\,\frac{d\!\vol_g}{r^4(x)} \right)^{1/2} \left( \int_{B(o,2^{2k})\setminus
B(o,2^k)}  \frac{G_o}{|\log G_o|^{2}} |d\log \bar G_o|^2\, d\!\vol_g \right)^{1/2}
\end{split}
\end{equation*}
The function $G_0$ tends to $0$ at infinity by Lemma 3.2 in \cite{kleinian},
hence we can again choose $R>0$ large enough so that $G_o\le 1/2$ outside the ball
$B(o,R/2)$. As a result, we have:
$$\int_{B(o,2^{2k})\setminus B(o,2^k)} G_o|\log G_o|^{2} 
\, \frac{d\!\vol_g}{r^4(x)} \ \le \ C\,k \ . $$
We now use another estimate due to R. Schoen and S.T. Yau which is valid in 
every dimension: they show at the end of page 60 in \cite{kleinian} that for 
any $\psi \in C^\infty_0(M)$ vanishing near $o$ we have :
$$\int_M \psi^2 |\nabla\log \bar G_o|^2  \, d\!\vol_g\ \le 
\ c\,\int_M \psi^2 |\nabla\log G_o|^2 \, d\!\vol_g 
\ + \ c\,\int_M |d\psi| ^2 \, d\!\vol_g$$
Let now $\phi$ a function with the following properties:
\begin{itemize}
\item[-] it has compact support in $B(o,2^{2k+1})\setminus B(o,2^{k-1})$;
\item[-] it is identically $1$ on $B(o,2^{2k})\setminus B(o,2^k)$;
\item[-] its gradient is bounded by $2^{-2k}\, c$ on 
$B(o,2^{2k+1})\setminus B(o,2^{2k})$;
\item[-] and its gradient is bounded by $2^{-k}\, c$ on 
$B(o,2^{k})\setminus B(o,2^{k-1})$.
\end{itemize}
This does exist and we can apply the last estimate to 
$\psi=(\phi\sqrt{G_o})/\log G_o$. We obtain :
\begin{equation*}
\begin{split} 
\int_{B(o,2^{2k})\setminus B(o,2^k)} & G_o|\log G_o|^{-2} |d\log \bar G_o|^2     
\, d\!\vol_g \\
& \le \ \ C\, \int_{B(o,2^{2k+1})\setminus
B(o,2^{k-1})} |d\phi|^2 G_o\, |\log G_o|^{-2}\, d\!\vol_g \\
& \ \ \ \ \ \ \ \ \ + C \,\int_{B(o,2^{2k+1})\setminus B(o,2^{k-1})} 
G_o|\, (\log G_o|^{-2} + |\log G_o|^{-4})\, |d\log  G_o|^2\, d\!\vol_g
\end{split}\end{equation*}
From Lemma \ref{lem:first} we know that the second term in the right-hand side
of this last estimate goes to zero when $k$ goes to infinity. 
For the first term, Cauchy-Schwartz inequality leads to :
\begin{equation*}
\begin{split}
& \int_{B(o,2^{2k+1})\setminus B(o,2^{k-1})} |d\phi|^2 \frac{G_o}{ |\log G_o|^{2} }
\, d\!\vol_g \\
& \ \ \ \ \ \leq \ \left( \int_{B(o,2^{2k+1})\setminus B(o,2^{k-1})} |d\phi|^4  
\, d\!\vol_g\right)^{1/2}\,\left(\int_{B(o,2^{2k+1})\setminus B(o,2^{k-1})} 
G_o^2 |\log G_o|^{-4} \, d\!\vol_g\right)^{1/2} 
\end{split}\end{equation*}
But our volume growth assumptions imply that
 $$\int_{B(o,2^{2k+1})\setminus B(o,2^{k-1})} |d\phi|^4 \, d\!\vol_g\le C$$ 
 where $C$ is a constant independent of $k$.
Using Lemma \ref{finite2}, we end with
$$\lim_{k\to\infty} \int_{B(o,2^{2k+1})\setminus B(o,2^{k-1})} 
G_o^2 |\log G_o|^{-4} \, d\!\vol_g  = 0 \ .$$
Hence we obtain again that $v$ is the constant function equal to $1$ everywhere. 
And this is what is needed  to resume the proof of \cite{kleinian}. \qed

\bigskip

\bibliographystyle{amsplain}
\providecommand{\bysame}{\leavevmode\hbox to3em{\hrulefill}\thinspace}

\bigskip

\end{document}